\input amstex
\input amsppt.sty
\magnification=\magstep1
\hsize=32truecc
\vsize=22.2truecm
\baselineskip=16truept
\NoBlackBoxes
\TagsOnRight \pageno=1 \nologo
\def\Z{\Bbb Z}
\def\N{\Bbb N}

\def\l{\left}
\def\r{\right}
\def\bg{\bigg}
\def\({\bg(}
\def\[{\bg\lfloor}
\def\){\bg)}
\def\]{\bg\rfloor}
\def\t{\text}
\def\f{\frac}

\def\p{\ (\roman{mod}\ p)}

\def\bi{\binom}
\def\eq{\equiv}

\def\ls{\leqslant}
\def\gs{\geqslant}
\def\mo{\roman{mod}}

\def\ve{\varepsilon}

\def\Proof{\noindent{\it Proof}}

\def\Remark{\medskip\noindent{\it  Remark}}

\def\Ack{\medskip\noindent {\bf Acknowledgment}}
\hbox {Nanjing Univ. J. Math. Biquarterly 40 (2023), no.\,1, 1--33.}
\bigskip
\topmatter
\title New congruences involving harmonic numbers\endtitle
\author Zhi-Wei Sun\endauthor
\leftheadtext{Zhi-Wei Sun}
\affil Department of Mathematics, Nanjing University\\
 Nanjing 210093, People's Republic of China
  \\  zwsun\@nju.edu.cn
  \\ {\tt http://maths.nju.edu.cn/$\sim$zwsun}
\endaffil
\abstract Let $p>3$ be a prime. For any $p$-adic integer $a$, we determine
$$\sum_{k=0}^{p-1}\bi{-a}k\bi{a-1}kH_k,\ \sum_{k=0}^{p-1}\bi{-a}k\bi{a-1}kH_k^{(2)},\ \sum_{k=0}^{p-1}\bi{-a}k\bi{a-1}k\f{H_k^{(2)}}{2k+1}$$
modulo $p^2$, where $H_k=\sum_{0<j\ls k}1/j$ and $H_k^{(2)}=\sum_{0<j\ls k}1/j^2$. In particular,
we show that
$$\gather\sum_{k=0}^{p-1}\bi{-a}k\bi{a-1}kH_k\eq(-1)^{\langle a\rangle_p}2\l(B_{p-1}(a)-B_{p-1}\r)\pmod p,
\\\sum_{k=0}^{p-1}\bi{-a}k\bi{a-1}kH_k^{(2)}\eq -E_{p-3}(a)\pmod p,
\\(2a-1)\sum_{k=0}^{p-1}\bi{-a}k\bi{a-1}k\f{H_k^{(2)}}{2k+1}\eq B_{p-2}(a)\pmod p,
\endgather$$
where $\langle a\rangle_p$ stands for the least nonnegative integer $r$ with $a\eq r\pmod{p}$, and
$B_n(x)$ and $E_n(x)$ denote the Bernoulli polynomial of degree $n$ and the Euler polynomial of degree $n$ respectively. We also pose some new conjectures on congruences.
\endabstract
\thanks 2020 {\it Mathematics Subject Classification}. \,Primary 11A07, 11B65;
Secondary  05A10, 11B68, 11B75.
\newline\indent The research was supported by the National Natural Science
Foundation of China (grant no. 11971222), and the initial version was posted to arXiv in 2014 with the ID
{\tt arXiv:1407.8465}.
\endthanks

\endtopmatter
\document

\heading{1. Introduction}\endheading

 A classical theorem of J. Wolstenholme [W] asserts that for any prime $p>3$ we have
$$\f12\bi{2p}p=\bi{2p-1}{p-1}\eq1\pmod{p^3},$$
which follows from the congruences
$$ H_{p-1}\eq0\pmod{p^2}\quad\t{and}\quad H_{p-1}^{(2)}\eq0\pmod p,$$
where
$$H_n:=\sum_{0<k\ls n}\f1k\quad\t{and}\quad H_n^{(2)}:=\sum_{0<k\ls n}\f1{k^2}\quad\t{for}\ n\in\N=\{0,1,2,\ldots\}.$$
Those $H_n\ (n\in\N)$ are the usual harmonic numbers, and those $H_n^{(2)}\ (n\in\N)$ are called second-order harmonic numbers.
For some congruences involving harmonic numbers, one may consult [Su12], [Su22] and [SZ].

In 2003, based on his analysis of the $p$-adic analogues of Gaussian hypergeometric series and the Calabi-Yau manifolds, F. Rodriguez-Villegas [RV] conjectured that for any prime $p>3$ we have
$$\gather\sum_{k=0}^{p-1}\f{\bi{2k}k^2}{16^k}\eq\l(\f{-1}p\r)\pmod{p^2},\ \ \sum_{k=0}^{p-1}\f{\bi{2k}k\bi{3k}k}{27^k}\eq\l(\f{p}3\r)\pmod{p^2},
\\\sum_{k=0}^{p-1}\f{\bi{4k}{2k}\bi{2k}k}{64^k}\eq\l(\f{-2}p\r)\pmod{p^2},\ \ \sum_{k=0}^{p-1}\f{\bi{6k}{3k}\bi{3k}k}{432^k}\eq\l(\f{-1}p\r)\pmod{p^2},
\endgather$$
where $(\f{\cdot}p)$ denotes the Legendre symbol. All the four congruences were proved by E. Mortenson [M1, M2] via the $p$-adic $\Gamma$-function and modular forms.
Z.-H. Sun [S1] presented elementary proofs of them, and V.J.W. Guo and J. Zeng [GZ] obtained a $q$-analogue of the first one.

Let $p>3$ be a prime. The author [Su11] showed that
$$\sum_{k=0}^{p-1}\f{\bi{2k}k^2}{16^k}\eq\l(\f{-1}p\r)-p^2E_{p-3}\pmod{p^3}\tag1.1$$
(see also [Su13] for a simpler proof), and conjectured that
$$\align\sum_{k=0}^{p-1}\f{\bi{2k}k\bi{3k}k}{27^k}\eq&\l(\f p3\r)-\f{p^2}3B_{p-2}\l(\f13\r)\pmod{p^3},\tag1.2
\\\sum_{k=0}^{p-1}\f{\bi{4k}{2k}\bi{2k}k}{64^k}\eq&\l(\f{-2}p\r)-\f3{16}p^2E_{p-3}\l(\f14\r)\pmod{p^3},\tag1.3
\\\sum_{k=0}^{p-1}\f{\bi{6k}{3k}\bi{3k}k}{432^k}\eq&\l(\f{-1}p\r)-\f{25}9p^2E_{p-3}\pmod{p^3},\tag1.4
\\\sum_{k=0}^{p-1}\f{\bi{2k}k\bi{3k}k}{(2k+1)27^k}\eq&\l(\f p3\r)-\f23p^2B_{p-2}\l(\f13\r)\pmod{p^3},\tag1.5
\\\sum_{k=0}^{p-1}\f{\bi{4k}{2k}\bi{2k}k}{(2k+1)64^k}\eq&\l(\f{-1}p\r)-3p^2E_{p-3}\pmod{p^3},\tag1.6
\\\sum_{k=0}^{p-1}\f{\bi{6k}{3k}\bi{3k}k}{(2k+1)432^k}\eq&\l(\f p3\r)\pmod{p^2},\tag1.7
\endalign$$
where $E_0,E_1,E_2,\ldots$ are the Euler numbers, and $E_n(x)$ denotes the Euler polynomial of degree $n$ given by
$$E_n(x)=\sum_{k=0}^n\bi nk\f{E_k}{2^k}\l(x-\f12\r)^{n-k},$$
and $B_n(x)$ stands for the Bernoulli polynomial of degree $n$ given by
$$B_n(x)=\sum_{k=0}^n\bi nk B_kx^{n-k}$$
with $B_0,B_1,B_2,\ldots$ the Bernoulli numbers. The conjectural congruences (1.2)-(1.7) were confirmed by Z.-H. Sun [S2, S3].

In this paper we mainly establish two new theorems involving harmonic numbers and second-order harmonic numbers.

For a prime $p$ and a $p$-adic integer $a$, we write $\langle a\rangle_p$ for the unique integer $r\in\{0,1,\ldots,p-1\}$ with $a\eq r\pmod p$,
and let $q_p(a)$ denote the Fermat quotient $(a^{p-1}-1)/p$ if $a\not\eq0\pmod p$.

\proclaim{Theorem 1.1} Let $p>3$ be a prime. For any $p$-adic integer $a$, we have
$$\aligned\sum_{k=0}^{p-1}\bi {-a}k\bi{a-1}kH_k\eq&(-1)^{\langle a\rangle_p-1}2\sum_{0<k<\langle a\rangle_p}\f1{a-k}
\\\eq&(-1)^{\langle a-1\rangle_p}\l(2H_{\langle a-1\rangle_p}+(a-\langle a\rangle_p)B_{p-2}(a)\r)\pmod{p^2}
\\\eq&(-1)^{\langle a\rangle_p}2(B_{p-1}(a)-B_{p-1})\pmod p.
\endaligned\tag1.8$$
\endproclaim
\Remark\ 1.1. Let $p>3$ be a prime and let $a$ be a $p$-adic integer.
Congruences involving the general sum $\sum_{k=0}^{p-1}\bi{a}k\bi{a+k}k/m^k=\sum_{k=0}^{p-1}\bi{a}k\bi{-1-a}k/(-m)^k$
with $m\not\eq0\pmod p$ first appeared in the author's paper [Su14]. Z.-H. Sun [S1, Corollary 2.1] determined
$\sum_{k=0}^{p-1}\bi ak\bi{-1-a}k$ modulo $p^2$ with the special cases $a=-1/2,-1/3,-1/4,-1/6$ first discovered by Rodriguez-Villegas [RV].
Besides Theorem 1.1, we are also able to show that
$$\sum_{k=1}^{p-1}\bi {-a}k\bi{a-1}k\f{H_k}k\eq(-1)^{\langle -a\rangle_p}E_{p-3}(a)\pmod p.$$
\medskip

 Let $p>3$  be a prime. As
$$\bi{p-1}k(-1)^k=\prod_{0<j\ls k}\l(1-\f pj\r)\eq1-pH_k\pmod{p^2}\quad\t{for all}\ k=0,1,2,\ldots,$$
combining Theorem 1.1 with [S1, Corollary 2.1], we obtain
$$\aligned&\sum_{k=0}^{p-1}\bi{p-1}k\bi {-a}k\bi{a-1}k(-1)^k
\\\eq&(-1)^{\langle-a\rangle_p}\l(1+2p\l(B_{p-1}(a)-B_{p-1}\r)\r)\pmod{p^2}
\endaligned\tag1.9$$
for any $p$-adic integer $a$.
For each $d=2,3,4,6$ and any $c\in\{1,\ldots,d\}$ with $(c,d)=1$, E. Lehmer [L] determined $B_{p-1}(c/d)-B_{p-1}$ modulo $p$ in terms of Fermat quotients.
For $d\in\{5,8,10,12\}$ and $c\in\{1,\ldots,d\}$ with $(c,d)=1$, A. Granville and the author [GS] determined $B_{p-1}(c/d)-B_{p-1}$ mod $p$ by showing that
$$\align
B_{p-1}\left(\frac{c}{5}\right) - B_{p-1}&\equiv \frac{5}{4}\(\left(
\frac{cp}{5}\right) \frac{1}{p} F_{p-(\frac{5}{p})} + q_p(5)\) \pmod p, \cr
B_{p-1}\left(\frac{c}{8}\right) - B_{p-1}&\equiv
\left(\frac{2}{cp}\right) \frac{2}{p} P_{p-(\frac{2}{p})} + 4q_p(2) \pmod p,
 \cr B_{p-1}\left(\frac{a}{10}\right) - B_{p-1}&\equiv \frac{15}{4}
\left(\frac{cp}{5}\right) \frac{1}{p} F_{p-(\frac{5}{p})} +
\frac{5}{4} q_p(5)+ 2q_p(2) \pmod p, \cr
B_{p-1}\left(\frac{c}{12}\right) - B_{p-1}&\equiv
\left(\frac{3}{c}\right) \frac{3}{p} S_{p-(\frac{3}{p})} + 3q_p(2)
 + \frac{3}{2}q_p(3) \pmod p,\endalign$$
where $(-)$ is the Jacobi symbol, and the Fibonacci sequence $(F_n)_{n\gs0}$, the Pell sequence $(P_n)_{n\gs0}$, and the sequence $(S_n)_{n\gs0}$
(cf. [Su02]) are defined as follows:
$$\align
&F_0 = 0, \ F_1 = 1, \ \ \t{and} \ F_{n+1} = F_{n} + F_{n-1} \ \ \ \t{for all} \ n=1,2,3,\ldots; \cr
&P_0 = 0, \ P_1 = 1, \ \ \t{and} \ P_{n+1} = 2P_{n} + P_{n-1} \ \  \t{for all} \ n=1,2,3,\ldots; \cr
&S_0 = 0, \ S_1 = 1, \ \ \t{and} \ S_{n+1} = 4S_{n} - S_{n-1} \ \ \t{for all} \ n=1,2,3,\ldots. \cr
\endalign$$

\proclaim{Corollary 1.1} Let $p>3$ be a prime. Then
$$\align\l(\f{-1}p\r)\sum_{k=0}^{p-1}\f{\bi{2k}k^2}{16^k}H_k\eq&-4q_p(2)+2p\,q_p(2)^2\pmod{p^2},\tag1.10
\\\l(\f p3\r)\sum_{k=0}^{p-1}\f{\bi{2k}k\bi{3k}k}{27^k}H_k\eq&-3q_p(3)+\f32p\,q_p(3)^2\pmod {p^2}\tag1.11
\\\l(\f{-2}p\r)\sum_{k=0}^{p-1}\f{\bi{4k}{2k}\bi{2k}k}{64^k}H_k\eq&-6q_p(2)+3p\,q_p(2)^2\pmod {p^2},\tag1.12
\endalign$$
and
$$\l(\f{-1}p\r)\sum_{k=0}^{p-1}\f{\bi{6k}{3k}\bi{3k}k}{432^k}H_k\eq-3q_p(3)-4q_p(2)+p\l(\f32q_p(3)^2+2q_p(2)^2\r)\,(\mo\ p^2).\tag1.13$$
\endproclaim

\proclaim{Theorem 1.2} Let $p>3$ be a prime, and let $a$ be a $p$-adic integer.

{\rm (i)} If $m$ is a positive integer with $a\eq m\pmod{p^2}$, then
$$\sum_{k=0}^{p-1}\bi{-a}k\bi{a-1}kH_k^{(2)}\eq 2\sum\Sb 0<k<m\\p\nmid k\endSb\f{(-1)^{m-k}}{k^2}\eq-E_{p^2-p-2}(a)\pmod{p^2}\tag1.14$$
and
$$\aligned&(2a-1)\sum_{k=0}^{p-1}\bi{-a}k\bi{a-1}k\f{H_k^{(2)}}{2k+1}
\\\eq&-2\sum\Sb 0<k<m\\p\nmid k\endSb\f1{k^2}\eq(2-2p)B_{p^2-p-1}(a)\pmod{p^2}.\endaligned\tag1.15$$

{\rm (ii)} We always have
$$\sum_{k=0}^{p-1}\bi{-a}k\bi{a-1}kH_k^{(2)}\eq -E_{p-3}(a)\pmod p\tag1.16$$
and
$$(2a-1)\sum_{k=0}^{p-1}\bi{-a}k\bi{a-1}k\f{H_k^{(2)}}{2k+1}\eq B_{p-2}(a)\pmod p.\tag1.17$$
\endproclaim
\Remark\ 1.2. Let $p>3$ be a prime. As $H_{(p-1)/2}^{(2)}\eq0\pmod p$, the number $H_k^{(2)}/(2k+1)$ is a $p$-adic integer for every $k=0,1,\ldots,p-1$.
For any $p$-adic integer $a$, Z.-H. Sun [S2] determined $\sum_{k=0}^{p-1}\bi {-a}k\bi{a-1}k$ and $\sum_{k=0}^{p-1}\bi {-a}k\bi{a-1}k\f1{2k+1}$
(with $a\not\eq1/2\pmod p$) modulo $p^3$. Combining this with Theorem 1.2(ii), we determine
$$\sum_{k=0}^{p-1}(-1)^k\bi{p-1}k\bi{p+k}k\bi {-a}k\bi{a-1}k$$
and $$\sum_{k=0}^{p-1}\f{(-1)^k}{2k+1}\bi{p-1}k\bi{p+k}k\bi {-a}k\bi{a-1}k$$
modulo $p^3$ since
$$(-1)^k\bi{p-1}k\bi{p+k}k=\prod_{0<j\ls k}\l(1-\f{p^2}{j^2}\r)\eq1-p^2H_k^{(2)}\pmod{p^4}\ \ \t{for}\ k\in\N.$$

\proclaim{Corollary 1.2} Let $p>3$ be a prime. Then
$$\aligned\sum_{k=0}^{p-1}\f{\bi{4k}{2k}\bi{2k}k}{64^k}H_k^{(2)}\eq&-E_{p^2-p-2}\l(\f14\r)\pmod {p^2}
\\\eq&-E_{p-3}\l(\f14\r)\pmod p.\endaligned\tag1.18$$
and
$$\aligned\sum_{k=0}^{p-1}\f{\bi{2k}k^2}{16^k}H_k^{(2)}\eq&\f14\sum_{k=0}^{p-1}\f{\bi{4k}{2k}\bi{2k}k}{(2k+1)64^k}H_k^{(2)}
\eq\f15\sum_{k=0}^{p-1}\f{\bi{6k}{3k}\bi{3k}k}{432^k}H_k^{(2)}
\\\eq&-4E_{p^2-p-2}\pmod{p^2}
\\\eq&-4E_{p-3}\pmod p.
\endaligned\tag1.19$$
We also have
$$\aligned \sum_{k=0}^{p-1}\f{\bi{2k}k\bi{3k}k}{27^k}H_k^{(2)}\eq&\f12\sum_{k=0}^{p-1}\f{\bi{2k}k\bi{3k}k}{(2k+1)27^k}H_k^{(2)}
\eq\f15\sum_{k=0}^{p-1}\f{\bi{6k}{3k}\bi{3k}k}{(2k+1)432^k}H_k^{(2)}\pmod{p^2}
\\\eq&(3p-3)B_{p^2-p-1}\l(\f13\r)\pmod{p^2}
\\\eq&-\f32B_{p-2}\l(\f13\r)\pmod p.
\endaligned\tag1.20$$
\endproclaim
\Remark\ 1.3. (i) The author [Su15] reported that $205129$ is the first odd prime $p$ with $B_{p-2}(1/3)\eq0\pmod p$.
We note that $1019$ is the first odd prime $p$ with $E_{p-3}(1/4)\eq0\pmod p$.

(ii) The author [S11, Conjecture 5.12(iii)] conjectured that for any prime $p>3$ we have
$$2\sum_{k=0}^{p-1}\f{\bi{2k}k\bi{3k}k}{27^k}-\sum_{k=0}^{p-1}\f{\bi{2k}k\bi{3k}k}{(2k+1)27^k}\eq\l(\f p3\r)\pmod {p^4}.\tag1.21$$

We are going to show Theorem 1.1 and Corollary 1.1 in the next section,
and prove Theorem 1.2 and Corollary 1.2 in Section 3. In Section 4, we pose some new conjectures on congruences.

\heading{2. Proofs of Theorem 1.1 and Corollary 1.1}\endheading

\proclaim{Lemma 2.1} For any positive integer $k$, we have the polynomial identity
$$\bi{-x}k\bi{x-1}k+\bi{x}k\bi{-x-1}k=2\l(\bi{x-1}k\bi{-x-1}k-\bi{x-1}{k-1}\bi{-x-1}{k-1}\r).\tag2.1$$
\endproclaim
\Proof. We may deduce (2.1) in view of [S1, p.\,310], but here we give a direct proof.
$$\align&\bi{-x}k\bi{x-1}k+\bi{x}k\bi{-x-1}k
\\=&\f{(-1)^k}{k!k!}\l((x-k)\cdots(x+k-1)+(x-k+1)\cdots(x+k)\r)
\\=&\f{(-1)^k}{k!k!}(x-k+1)\cdots(x+k-1)(x-k+x+k)
\\=&2\f{(-1)^k}{k!k!}\cdot\f{(x-k+1)\cdots(x+k-1)}{x}\l((x-k)(x+k)+k^2\r)
\\=&2\f{(-1)^k}{k!k!}\l(\f{(x-k)\cdots(x+k)}{x}+k^2\f{(x-k+1)\cdots(x+k-1)}{x}\r)
\\=&2\l(\bi{x-1}k\bi{-x-1}k-\bi{x-1}{k-1}\bi{-x-1}{k-1}\r).
\endalign$$
This completes the proof. \qed

\proclaim{Lemma 2.2} For any positive integer $n$, we have
$$\f1n\sum_{k=1}^n(k^2-kx^2)\bi xk\bi{-x}k=(n^2-x^2)\bi xn\bi{-x}n.\tag2.2$$
\endproclaim
\Proof. It is easy to verify (2.2) for $n=1$.

Now assume that (2.2) holds for a fixed positive integer $n$. Then
$$\align &\sum_{k=1}^{n+1}(k^2-kx^2)\bi xk\bi{-x}k
\\=&n(n^2-x^2)\bi xn\bi{-x}n+((n+1)^2-(n+1)x^2)\bi x{n+1}\bi{-x}{n+1}
\\=&\l(n(n+1)^2+(n+1)^2-(n+1)x^2\r)\bi x{n+1}\bi{-x}{n+1}
\\=&(n+1)\l((n+1)^2-x^2\r)\bi x{n+1}\bi{-x}{n+1}.
\endalign$$
This concludes the induction proof. \qed

\medskip
\noindent{\it Proof of Theorem} 1.1. Define
$$P_n(x):=\sum_{k=0}^n\bi{-x}k\bi{x-1}kH_k\quad\t{for}\ n=0,1,2,\ldots.$$
With the help of (2.1), we have
$$\align &P_n(x)+P_n(x+1)
\\=&\sum_{k=1}^n\(\bi{-x}k\bi{x-1}k+\bi xk\bi{-x-1}k\)H_k
\\=&2\sum_{k=1}^n\(\bi{x-1}k\bi{-x-1}kH_k-\bi{x-1}{k-1}\bi{-x-1}{k-1}\l(H_{k-1}+\f1k\r)\)
\\=&2\bi{x-1}n\bi{-x-1}nH_n+\f 2{x^2}\sum_{k=1}^n k\bi xk\bi{-x}k.
\endalign$$
Recall that $H_{p-1}\eq0\pmod{p^2}$. Thus, for any $p$-adic integer $x\not\eq0\pmod{p}$, we have
$$P_{p-1}(x)+P_{p-1}(x+1)\eq\f 2{x^2}\sum_{k=1}^{p-1} k\bi xk\bi{-x}k\pmod{p^2}.\tag2.3$$
If $x\eq0\pmod p$, then
$$\align P_{p-1}(x)=&\sum_{k=1}^{p-1}\f{-x}k\bi{-x-1}{k-1}\bi{x-1}kH_k
\\\eq&-x\sum_{k=1}^{p-1}\f1k\bi{-1}{k-1}\bi{-1}kH_k=x\sum_{k=1}^{p-1}\l(\f1{k^2}+\f{H_{k-1}}k\r)
\\\eq&x\sum_{0<j<k<p}\f1{jk}=\f x2\l(H_{p-1}^2-H_{p-1}^{(2)}\r)\eq0\pmod{p^2}
\endalign$$
and also $P_{p-1}(x+1)=P_{p-1}(-x)\eq0\pmod{p^2}$.

In light of (2.2), for any positive integer $n$ we have
$$\align &x^2\sum_{k=1}^nk\bi xk\bi{-x}k+n(n^2-x^2)\bi xn\bi{-x}n
\\=&\sum_{k=1}^nk^2\bi xk\bi{-x}k=-x^2\sum_{k=1}^n\bi{x-1}{k-1}\bi{-x-1}{k-1}
\\=&x^2\bi{x-1}n\bi{-x-1}n-x^2\sum_{k=0}^n\bi{x-1}k\bi{-x-1}k
\\=&-(n+1)^2\bi x{n+1}\bi{-x}{n+1}-x^2\sum_{k=0}^n\bi{x-1}k\bi{-x-1}k.
\endalign$$
Let $x$ be any $p$-adic integer with $x\not\eq0\pmod p$. Clearly,
$$\bi{x}{p-1}\bi{-x}{p-1}=\f{\prod_{r=0}^{p-2}(r^2-x^2)}{((p-1)!)^2}\eq0\pmod p$$
and hence
$$\align&\l((p-1)^2-x^2\r)\bi x{p-1}\bi{-x}{p-1}
\\\eq&(1-x^2)\bi x{p-1}\bi{-x}{p-1}=p^2\bi{x+1}p\bi{1-x}p\eq0\pmod{p^2}.\endalign$$
Thus
$$\align \sum_{k=1}^{p-1}k\bi xk\bi{-x}k
\eq&-\sum_{k=0}^{p-1}\bi{x-1}k\bi{-x-1}k
\\\eq&-\bi{p-2}{\langle x-1\rangle_p}\l(1+(-2-(p-2))H_{p-2}\r)
\\&+\bi{p-2}{\langle x-1\rangle_p}(x-1-\langle x-1\rangle_p)H_{\langle x-1\rangle_p}
\\&+\bi{p-2}{\langle x-1\rangle_p}(-x-1-\langle-x-1\rangle_p)H_{\langle-x-1\rangle_p}\pmod{p^2},
\endalign$$
with the help of [S1, Theorem 4.1].
Note that
$$\gather 1-pH_{p-2}=1+\f p{p-1}-pH_{p-1}\eq1-p\pmod{p^2},
\\H_{\langle-x-1\rangle_p}=H_{p-1-\langle x\rangle_p}=H_{p-1}-\sum_{k=1}^{\langle x\rangle_p}\f1{p-k}\eq H_{\langle x\rangle_p}\pmod p
\endgather$$
and
$$\align(p-1)\bi{p-2}{\langle x-1\rangle_p}=&\langle x\rangle_p\bi{p-1}{\langle x\rangle_p}
=(-1)^{\langle x\rangle_p}\langle x\rangle_p\prod_{k=1}^{\langle x\rangle_p}\l(1-\f pk\r)
\\\eq&(-1)^{\langle x\rangle_p}\langle x\rangle_p\l(1-pH_{\langle x\rangle_p}\r)\pmod{p^2}
\\\eq&(-1)^{\langle x\rangle_p}x\pmod p.
\endalign$$
Combining this with (2.3), we get
$$\align&\f{x^2}2\l(P_{p-1}(x)+P_{p-1}(x+1)\r)
\\\eq&\sum_{k=1}^{p-1}k\bi xk\bi{-x}k
\\\eq&(-1)^{\langle x\rangle_p}\langle x\rangle_p\l(1-pH_{\langle x\rangle_p}\r)-x(-1)^{\langle x\rangle_p}(x-\langle x\rangle_p)H_{\langle x\rangle_p-1}
\\&+x(-1)^{\langle x\rangle_p}(x+1+p-1-\langle x\rangle_p)H_{\langle x\rangle_p}
\\\eq&(-1)^{\langle x\rangle_p}x\pmod{p^2}.
\endalign$$

By the above, for any $p$-adic integer $x$, we have
$$P_{p-1}(x)+P_{p-1}(x+1)
\eq\cases(-1)^{\langle x\rangle_p}2/x\pmod{p^2}&\t{if}\ x\not\eq0\pmod p,
\\0\pmod{p^2}&\t{otherwise}.\endcases\tag2.4$$
Therefore
$$\align -P_{p-1}(a)\eq&(-1)^{\langle a\rangle_p}P_{p-1}(a-\langle a\rangle_p)-P_{p-1}(a)
\\=&\sum_{0<k\ls\langle a\rangle_p}\l((-1)^kP_{p-1}(a-k)-(-1)^{k-1}P_{p-1}(a-k+1)\r)
\\\eq&\sum_{0<k<\langle a\rangle_p}(-1)^k(-1)^{\langle a-k\rangle_p}\f2{a-k}=2(-1)^{\langle a\rangle_p}\sum_{0<k<\langle a\rangle_p}\f1{a-k}
\\=&2(-1)^{\langle a\rangle_p}\sum_{0<k<\langle a\rangle_p}\(\f1{\langle a\rangle_p-k}
+\f{\langle a\rangle_p-k-(a-k)}{(a-k)(\langle a\rangle_p-k)}\)
\\\eq&2(-1)^{\langle a\rangle_p}\l(H_{\langle a-1\rangle_p}+(\langle a\rangle_p-a)H_{\langle a-1\rangle_p}^{(2)}\r)\pmod{p^2}
\endalign$$
and hence
$$\align P_{p-1}(a)\eq&2(-1)^{\langle a\rangle_p-1}\sum_{0\ls k<\langle a\rangle_p}k^{p-2}
\\=&\f{2(-1)^{\langle a\rangle_p-1}}{p-1}\sum_{0\ls k<\langle a\rangle_p}\l(B_{p-1}(k+1)-B_{p-1}(k)\r)
\\\eq&2(-1)^{\langle a\rangle_p}\l(B_{p-1}(\langle a\rangle_p)-B_{p-1}\r)
\\\eq&2(-1)^{\langle a\rangle_p}\l(B_{p-1}(a)-B_{p-1}\r) \pmod p.
\endalign$$
Note also that
$$\align H_{\langle a-1\rangle_p}^{(2)}\eq&\sum_{0\ls k<\langle a\rangle_p}k^{p-3}
=\sum_{0\ls k<\langle a\rangle_p}\f{B_{p-2}(k+1)-B_{p-2}(k)}{p-2}
\\=&\f{B_{p-2}(\langle a\rangle_p)-B_{p-2}}{p-2}\eq-\f12B_{p-2}(a)\pmod p.
\endalign$$
So we have the desired (1.8). \qed

The following lemma was first deduced by E. Lehmer [L].

\proclaim{Lemma 2.3} Let $p>3$ be a prime. Then
$$\align\sum_{k=1}^{\lfloor p/2\rfloor}\f1{p-2k}\eq& q_p(2)-\f p2q_p(2)^2\pmod{p^2},\tag2.5
\\\sum_{k=1}^{\lfloor p/3\rfloor}\f1{p-3k}\eq& \f{q_p(3)}2-\f p4q_p(3)^2\pmod{p^2},\tag2.6
\\\sum_{k=1}^{\lfloor p/4\rfloor}\f1{p-4k}\eq& \f{3}4q_p(2)-\f 38p\,q_p(2)^2\pmod{p^2}.\tag2.7
\endalign$$
If $p>5$, then
$$\sum_{k=1}^{\lfloor p/6\rfloor}\f1{p-6k}\eq \f{q_p(3)}4+\f{q_p(2)}3-p\l(\f{q_p(3)^2}8+\f{q_p(2)^2}6\r)\pmod{p^2}.\tag2.8$$
\endproclaim

\medskip
\noindent{\it Proof of Corollary} 1.1. It is well known that for any $k\in\N$ we have
$$\align\bi{-1/2}k\bi{-1/2}k=&\l(\f{\bi{2k}k}{(-4)^k}\r)^2=\f{\bi{2k}k^2}{16^k},
\\\bi{-1/3}k\bi{-2/3}k=&\f{\bi{2k}k\bi{3k}k}{27^k},
\\\bi{-1/4}k\bi{-3/4}k=&\f{\bi{4k}{2k}\bi{2k}k}{64^k},
\\\bi{-1/6}k\bi{-5/6}k=&\f{\bi{6k}{3k}\bi{3k}k}{432^k}.
\endalign$$
Applying the first congruence in (1.8) with $a=1/2$ and Lehmer's congruence (2.5), we obtain
$$\align\sum_{k=0}^{p-1}\f{\bi{2k}k^2}{16^k}H_k\eq&(-1)^{(p-1)/2}2\sum_{k=1}^{(p-1)/2}\f1{1/2-k}
\\=&-4\l(\f{-1}p\r)\sum\Sb j=1\\4\nmid j\endSb^{p-1}\f1j=-4\l(\f{-1}p\r)\sum_{k=1}^{\lfloor p/2\rfloor}\f1{p-2k}
\\\eq&-4\l(\f{-1}p\r)\l(q_p(2)-\f p2q_p(2)^2\r)\pmod{p^2}.
\endalign$$
This proves (1.10). Choose $r\in\{1,2\}$ with $r\eq-p\pmod 3$. Then $\langle r/3\rangle_p=(p+r)/3$. By the first congruence in (1.8) with $a=r/3$ and Lehmer's congruence (2.6), we have
$$\align\sum_{k=0}^{p-1}\f{\bi{2k}k\bi{3k}k}{27^k}H_k\eq&(-1)^{(p+r)/3-1}2\sum_{0<k<(p+r)/3}\f1{r/3-k}
\\=&-6(-1)^r\sum\Sb j=1\\3\mid j+r\endSb^{p-1}\f1j
=-6\l(\f{p}3\r)\sum_{k=1}^{\lfloor p/3\rfloor}\f1{p-3k}
\\\eq&-6\l(\f{p}3\r)\l(\f{q_p(3)}2-\f p4q_p(3)^2\r)\pmod{p^2}.
\endalign$$
This proves (1.11).
Choose $s\in\{1,3\}$ with $s\eq-p\pmod 4$. Then $\langle s/4\rangle_p=(p+s)/4$. By the first congruence in (1.8) with $a=s/4$ and Lehmer's congruence (2.7), we have
$$\align\sum_{k=0}^{p-1}\f{\bi{4k}{2k}\bi{2k}k}{64^k}H_k\eq&(-1)^{(p+s)/4-1}2\sum_{0<k<(p+s)/4}\f1{s/4-k}
\\=&8(-1)^{(p+s)/4}\sum\Sb j=1\\4\mid j+s\endSb^{p-1}\f1j
=-8\l(\f{-2}p\r)\sum_{k=1}^{\lfloor p/4\rfloor}\f1{p-4k}
\\\eq&-8\l(\f{-2}p\r)\l(\f34q_p(2)-\f 38p\,q_p(2)^2\r)\pmod{p^2}.
\endalign$$
This proves (1.12). Choose $t\in\{1,5\}$ with $t\eq-p\pmod 6$. Then $\langle t/6\rangle_p=(p+t)/6$.
Provided $p>5$, by the first congruence in (1.8) with $a=t/6$ and Lehmer's congruence (2.8), we have
$$\align\sum_{k=0}^{p-1}\f{\bi{6k}{3k}\bi{3k}k}{432^k}H_k\eq&(-1)^{(p+t)/6-1}2\sum_{0<k<(p+t)/6}\f1{t/6-k}
\\=&12(-1)^{(p+t)/6}\sum\Sb j=1\\6\mid j+t\endSb^{p-1}\f1j
=-12\l(\f{-1}p\r)\sum_{k=1}^{\lfloor p/6\rfloor}\f1{p-6k}
\\\eq&-12\l(\f{-1}p\r)\l(\f{q_p(3)}4+\f{q_p(2)}3-p\l(\f{q_p(3)^2}8+\f{q_p(2)^2}6\r)\r)\pmod{p^2}.
\endalign$$
This proves (1.13). (Note that (1.13) for $p=5$ can be verified directly.)  We are done. \qed

\heading{3. Proofs of Theorem 1.2 and Corollary 1.2}\endheading

For any $n\in\N$, we define
$$W_n(x):=\sum_{k=0}^{n}\bi{-x}k\bi{x-1}kH_k^{(2)}\quad\t{and}\ w_n(x):=\sum_{k=0}^n\bi{-x}k\bi{x-1}k\f{H_k^{(2)}}{2k+1}.\tag3.1$$

\proclaim{Lemma 3.1} For any $n\in\N$£¬ we have
$$W_n(x)+W_n(x+1)=2\bi{x-1}n\bi{-x-1}n\l(H_n^{(2)}+\f1{x^2}\r)-\f 2{x^2}\tag3.2$$
and
$$(2x+1)w_n(x+1)-(2x-1)w_n(x)=2\bi{x-1}n\bi{-x-1}n\l(H_n^{(2)}+\f1{x^2}\r)-\f2{x^2}.\tag3.3$$
\endproclaim
\Proof. For any positive integer $k$, there is the polynomial identity
$$\aligned&(2x+1)\bi xk\bi{-x-1}k-(2x-1)\bi{-x}k\bi{x-1}k
\\=&2(2k+1)\l(\bi{x-1}k\bi{-x-1}k-\bi{x-1}{k-1}\bi{-x-1}{k-1}\r).
\endaligned\tag3.4$$
In fact,
$$\align&(2x+1)\bi xk\bi{-x-1}k-(2x-1)\bi{-x}k\bi{x-1}k
\\=&\f{(-1)^k}{k!k!}\l((2x+1)(x-k+1)\cdots(x+k)-(2x-1)(x-k)\cdots(x+k-1)\r)
\\=&\f{(-1)^k}{k!k!}(x-k+1)\cdots(x+k-1)\l((2x+1)(x+k)-(2x-1)(x-k)\r)
\\=&(-1)^k\f{2(2k+1)}{k!k!}\cdot\f{(x-k+1)\cdots(x+k-1)}x\l((x-k)(x+k)+k^2\r)
\\=&2(2k+1)\l(\bi{x-1}k\bi{-x-1}k-\bi{x-1}{k-1}\bi{-x-1}{k-1}\r).
\endalign$$

In light of (2.1) and (3.4),
$$\align &W_n(x)+W_n(x+1)
\\=&(2x+1)w_n(x)-(2x-1)w_n(x+1)
\\=&2\sum_{k=1}^{n}\(\bi{x-1}k\bi{-x-1}kH_k^{(2)}-\bi{x-1}{k-1}\bi{-x-1}{k-1}\l(H_{k-1}^{(2)}+\f1{k^2}\r)\)
\\=&2\bi{x-1}{n}\bi{-x-1}{n}H_{n}^{(2)}-2\sum_{k=1}^{n}\f1{k^2}\bi{x-1}{k-1}\bi{-x-1}{k-1}
\\=&2\bi{x-1}{n}\bi{-x-1}{n}H_{n}^{(2)}+\f 2{x^2}\sum_{k=1}^{n}\bi{x}{k}\bi{-x}{k}.
\endalign$$
Combining this with the identity
$$\sum_{k=0}^n\bi xk\bi{-x}k=\bi{x-1}n\bi{-x-1}n\tag3.5$$
we immediately obtain (3.2) and (3.3).
Note that the polynomial identity (3.5) holds if and only if it is valid for all $x=-n,-n-1\ldots$.
For each $x=-n,-n-1,\ldots$, the identity (3.5) has the equivalent form
$$\sum_{k=0}^n\bi xk\bi{-x}{-x-k}=\bi{x-1}n\bi{-x-1}n$$
which is a special case of Andersen's identity
$$m\sum_{k=0}^n\bi xk\bi{-x}{m-k}=(m-n)\bi{x-1}n\bi{-x}{m-n}\ \ (m\gs n\gs0)\tag3.6$$
(cf. (3.14) of [G, p.\,23]). This concludes the proof.
\qed

\proclaim{Lemma 3.2} Let $p$ be any prime, and let $x$ be a nonzero $p$-adic integer.
Then
$$\bi{x-1}{p-1}\bi{-x-1}{p-1}\l(H_{p-1}^{(2)}+\f1{x^2}\r)-\f1{x^2}\eq\cases-1/x^2\pmod{p^2}&\t{if}\ x\not\eq0\pmod p,
\\0\pmod {p^2}&\t{otherwise}.\endcases\tag3.7$$
\endproclaim
\Proof. If $x\not\eq0\pmod p$, then
$$\bi{x-1}{p-1}\bi{-x-1}{p-1}=\f{p^2}{-x^2}\bi{x}p\bi{-x}p\eq0\pmod{p^2}$$
and hence (3.7) holds.

Below we assume $x\eq0\pmod p$. Write $x=p^nx_0$, where $n$ is a positive integer and $x_0$ is a $p$-adic integer with $x_0\not\eq0\pmod p$.
Clearly,
$$\align\bi{x-1}{p-1}\bi{-x-1}{p-1}=&\prod_{k=1}^{p-1}\l(\f{p^nx_0-k}k\cdot\f{-p^nx_0-k}k\r)=\prod_{k=1}^{p-1}\l(1-\f{p^{2n}x_0^2}{k^2}\r)
\\\eq&1-\sum_{k=1}^{p-1}\f{p^{2n}x_0^2}{k^2}=1-x^2H_{p-1}^{(2)}\pmod{p^{4n}}.
\endalign$$
and hence
$$\f{\bi{x-1}{p-1}\bi{-x-1}{p-1}-1}{x^2}\eq-H_{p-1}^{(2)}\pmod{p^{2n}}.$$
Therefore,
$$\bi{x-1}{p-1}\bi{-x-1}{p-1}\l(H_{p-1}^{(2)}+\f1{x^2}\r)-\f1{x^2}\eq H_{p-1}^{(2)}+\f{\bi{x-1}{p-1}\bi{-x-1}{p-1}-1}{x^2}\eq0\ (\mo\ p^2).$$
This completes the proof. \qed

\proclaim{Lemma 3.3} For any prime $p>3$, we have
$$\sum\Sb k=1\\p\nmid k\endSb^{(p^2-1)/2}\f1{k^2}\eq\sum\Sb k=1\\ p\nmid k\endSb^{p^2-1}\f1{k^2}\eq0\pmod{p^2}.\tag3.8$$
\endproclaim
\Proof. Since $\{2j:\ 0<j<p^2\ \&\ p\nmid j\}$ is a reduced system of residues modulo $p^2$, we have
$$\sum^{p^2-1}\Sb j=1\\p\nmid j\endSb\f1{(2j)^2}\eq\sum^{p^2-1}\Sb k=1\\p\nmid k\endSb\f1{k^2}\pmod{p^2}
\ \ \ \t{and hence}\ \ \ \sum^{p^2-1}\Sb k=1\\p\nmid k\endSb\f1{k^2}\eq0\pmod{p^2}.$$
Note also that
$$2\sum^{(p^2-1)/2}\Sb k=1\\p\nmid k\endSb\f1{k^2}\eq\sum^{(p^2-1)/2}\Sb k=1\\p\nmid k\endSb\l(\f1{k^2}+\f1{(p^2-k)^2}\r)=\sum^{p^2-1}\Sb k=1\\p\nmid k\endSb\f1{k^2}\pmod{p^2}.$$
Therefore (3.8) holds. \qed

\medskip
\noindent{\it Proof of Theorem} 1.2. (i) Let $m$ be any positive integer with $a\eq m\pmod{p^2}$. Obviously,
$$W_{p-1}(a)\eq W_{p-1}(m)\pmod{p^2}\quad\t{and}\quad w_{p-1}(a)\eq w_{p-1}(m)\pmod{p^2}.$$
In light of Lemmas 3.1 and 3.2,
$$\align&(-1)^mW_{p-1}(m)-(-1)^1W_{p-1}(1)
\\=&\sum_{0<k<m}\l((-1)^{k+1}W_{p-1}(k+1)-(-1)^kW_{p-1}(k)\r)
\\\eq&\sum\Sb 0<k<m\\p\nmid k\endSb(-1)^{k+1}\f{-2}{k^2}\pmod{p^2}
\endalign$$
and
$$\align&(2m-1)w_{p-1}(m)-(2\times1-1)w_{p-1}(1)
\\=&\sum_{0<k<m}\l((2k+1)w_{p-1}(k+1)-(2k-1)w_{p-1}(k)\r)
\\\eq&\sum\Sb 0<k<m\\p\nmid k\endSb\f{-2}{k^2}\pmod{p^2}.
\endalign$$
Note that $W_{p-1}(1)=0=w_{p-1}(1)$. So we have the first congruence in (1.14) as well as the first congruence in (1.15).

It is well-known that $E_n(x)+E_n(x+1)=2x^n$ and $E_{2n+2}(0)=0$ for all $n\in\N$. Let $\varphi$ be Euler's totient function. Then
$$\align2\sum\Sb 0<k<m\\p\nmid k\endSb\f{(-1)^{m-k}}{k^2}\eq&2\sum_{k=0}^{m-1}(-1)^{m-k}k^{\varphi(p^2)-2}
\\=&(-1)^m\sum_{k=0}^{m-1}(-1)^k\l(E_{\varphi(p^2)-2}(k)+E_{\varphi(p^2)-2}(k+1)\r)
\\=&(-1)^m\sum_{k=0}^{m-1}\l((-1)^kE_{\varphi(p^2)-2}(k)-(-1)^{k+1}E_{\varphi(p^2)-2}(k+1)\r)
\\=&(-1)^m\l(E_{\varphi(p^2)-2}(0)-(-1)^mE_{\varphi(p^2)-2}(m)\r)=-E_{\varphi(p^2)-2}(m)
\\\eq&-E_{\varphi(p^2)-2}(a)\pmod{p^2}.
\endalign$$
This proves the second congruence in (1.14).

To complete the proof of (1.15), we only need to show that
$$\sum\Sb 0<k<m\\p\nmid k\endSb\f1{k^2}\eq(p-1)B_{\varphi(p^2)-1}(a)\pmod{p^2}.$$
Suppose that $a\eq m'=m+p^2q\pmod{p^3}$ with $q\in\{1,\ldots,p\}$.
 Then
 $$\align&\sum\Sb 0<k<m'\\p\nmid k\endSb\f1{k^2}-\sum\Sb 0<k<m\\p\nmid k\endSb\f1{k^2}
 \\=&\sum\Sb 0<k<p^2q\\p\nmid m+k\endSb\f1{(m+k)^2}=\sum\Sb r=0\\p\nmid m+r\endSb^{p^2-1}\sum_{s=0}^{q-1}
 \f1{(m+r+p^2s)^2}
 \\\eq&q\sum^{p^2-1}\Sb r=0\\p\nmid m+r\endSb\f1{(m+r)^2}\eq q\sum^{p^2-1}\Sb k=1\\p\nmid k\endSb\f1{k^2}
 \eq0\pmod{p^2}
 \endalign$$
 by Lemma 3.3. On the other hand, we have
$$\sum\Sb 0<k<m'\\p\nmid k\endSb\f1{k^2}\eq(p-1)B_{\varphi(p^2)-1}(a)\pmod{p^2}.$$
In fact, as $B_n(x+1)-B_n(x)=nx^{n-1}$ and $B_{2n+3}=0$ for all $n\in\N$, and $pB_n$ is a $p$-adic integer for any $n\in\N$, we have
$$\align \sum\Sb 0<k<m'\\p\nmid k\endSb\f1{k^2}\eq&\sum_{k=0}^{m'-1}k^{\varphi(p^2)-2}=\sum_{k=0}^{m'-1}\f{B_{\varphi(p^2)-1}(k+1)-B_{\varphi(p^2)-1}(k)}{\varphi(p^2)-1}
\\=&\f{B_{\varphi(p^2)-1}(m')-B_{\varphi(p^2)-1}}{\varphi(p^2)-1}=\f{B_{\varphi(p^2)-1}(m')}{p^2-p-1}
\\\eq& (p-1)B_{\varphi(p^2)-1}(a)+(p-1)\sum_{k=0}^{\varphi(p^2)-1}\bi{\varphi(p^2)-1}k(pB_{\varphi(p^2)-1-k})\f{(m')^k-a^k}p
\\\eq&(p-1)B_{\varphi(p^2)-1}(a)\pmod{p^2}.
\endalign$$
Therefore (1.15) also holds.

(ii) Choose $m\in\{1,2,\ldots,p^2\}$ such that $a\eq m\pmod{p^2}$. Write $m=ps+r$ with $s\in\{0,\ldots,p-1\}$ and $r\in\{1,\ldots,p\}$.
Then, for any $\ve=\pm1$ we have
$$\align\sum\Sb0<k<m\\p\nmid k\endSb\f{\ve^k}{k^2}=&\sum_{k=1}^s\sum_{t=1}^{p-1}\f{\ve^{pk-t}}{(pk-t)^2}+\sum_{0<t<r}\f{\ve^{ps+t}}{(ps+t)^2}
\\\eq&\sum_{k=1}^s\ve^k\sum_{t=1}^{p-1}\f{\ve^t}{t^2}+\ve^s\sum_{0<t<r}\f{\ve^t}{t^2}
\\\eq&\ve^s\sum_{t=0}^{r-1}\ve^tt^{p-3}
\pmod{p}
\endalign$$
since $H_{p-1}^{(2)}\eq0\pmod p$ and
$$\sum_{t=1}^{p-1}\f{(-1)^t}{t^2}=\sum_{t=1}^{(p-1)/2}\l(\f{(-1)^t}{t^2}+\f{(-1)^{p-t}}{(p-t)^2}\r)\eq0\pmod p.$$
Thus, we deduce from (1.14) that
$$\align W_{p-1}(a)\eq&2(-1)^m\sum\Sb 0<k<m\\p\nmid k\endSb\f{(-1)^k}{k^2}\eq(-1)^r\sum_{t=0}^{r-1}(-1)^t(2t^{p-3})
\\=&(-1)^r\sum_{t=0}^{p-1}\l((-1)^tE_{p-3}(t)-(-1)^{t+1}E_{p-3}(t+1)\r)
\\=&(-1)^r\l(E_{p-3}(0)-(-1)^rE_{p-3}\r)=-E_{p-3}(r)
\\\eq&-E_{p-3}(a)\pmod p.
\endalign$$
Similarly, from (1.15) we obtain
$$\align (2a-1)w_{p-1}(a)\eq&-2\sum\Sb 0<k<m\\p\nmid k\endSb\f{1}{k^2}\eq\sum_{t=0}^{r-1}((p-2)t^{p-3})
\\=&\sum_{t=0}^{r-1}\l(B_{p-2}(t+1)-B_{p-2}(t)\r)=B_{p-2}(r)
\\\eq&B_{p-2}(a)\pmod p.
\endalign$$
Therefore both (1.16) and (1.17) are valid.

By the above, we have completed the proof of Theorem 1.2. \qed

\medskip
\noindent{\it Proof of Corollary} 1.2. Applying (1.14) and (1.16) with $a=1/4$ we immediately get (1.18).

For every $a=1,2,\ldots$, we obviously have $$E_{\varphi(p^a)-2}\l(\f12\r)=\f{E_{\varphi(p^a)-2}}{2^{\varphi(p^a)-2}}\eq4E_{\varphi(p^a)-2}\pmod{p^a}.$$
(1.14) and (1.16) with $a=1/2$ yield
$$\align\sum_{k=0}^{p-1}\f{\bi{2k}k^2}{16^k}H_k^{(2)}\eq&-4E_{p^2-p-2}\pmod{p^2}
\\\eq&-4E_{p-3}\pmod p.\endalign$$
Clearly $1/2\eq(p^2+1)/2\pmod{p^2}$ and $1/4\eq(3p^2+1)/4\pmod{p^2}$.
Applying Theorem 1.2(i) with $a=1/2,1/4$, we get
$$\sum_{k=0}^{p-1}\f{\bi{2k}k^2}{16^k}H_k^{(2)}\eq-2\sum^{(p^2-1)/2}\Sb k=1\\p\nmid k\endSb\f{(-1)^k}{k^2}\pmod{p^2}\tag3.9$$
and
$$-\f12\sum_{k=0}^{p-1}\f{\bi{4k}{2k}\bi{2k}k}{(2k+1)64^k}H_k^{(2)}\eq-2\sum^{3(p^2-1)/4}\Sb k=1\\p\nmid k\endSb\f1{k^2}\pmod{p^2}.\tag3.10$$
In view of (3.8),
$$\align\sum^{\f34(p^2-1)}\Sb k=1\\p\nmid k\endSb\f1{k^2}=&\sum^{p^2-1}\Sb k=1\\p\nmid k\endSb\f1{k^2}-\sum^{(p^2-1)/4}\Sb k=1\\p\nmid k\endSb\f1{(p^2-k)^2}
\\\eq&-2\sum^{(p^2-1)/4}\Sb k=1\\p\nmid k\endSb\f 2{(2k)^2}=-2\sum^{(p^2-1)/2}\Sb k=1\\p\nmid k\endSb\f{1+(-1)^k}{k^2}
\\\eq&-2\sum^{(p^2-1)/2}\Sb k=1\\p\nmid k\endSb\f{(-1)^k}{k^2}\pmod{p^2}.
\endalign$$
Combining this with (3.9) and (3.10) we obtain the first congruence in (1.19).
It is known that
$$E_n\l(\f16\r)=2^{-n-1}(1+3^{-n})E_n\quad\t{for all}\ n=0,2,4,6,\ldots$$
(see, e.g., G. J. Fox [F]). Thus, applying (1.14) with $a=1/6$ we get
$$\align\sum_{k=0}^{p-1}\f{\bi{6k}{3k}\bi{3k}k}{432^k}H_k^{(2)}\eq&-E_{p^2-p-2}\l(\f16\r)=-2^{-\varphi(p^2)+1}\l(1+3^{-\varphi(p^2)+2}\r)E_{p^2-p-2}
\\\eq&-20E_{p^2-p-2}\pmod{p^2}
\endalign$$
and hence
$$\f15\sum_{k=0}^{p-1}\f{\bi{6k}{3k}\bi{3k}k}{432^k}H_k^{(2)}\eq-4E_{p^2-p-2}\eq\sum_{k=0}^{p-1}\f{\bi{2k}k^2}{16^k}H_k^{(2)}\pmod{p^2}.$$
(The first congruence in the last formula can be verified directly for $p=5$.)
This concludes the proof of (1.19).

 By (1.15) and (1.17) in the case $a=1/3$, we have
$$\align\sum_{k=0}^{p-1}\f{\bi{2k}k\bi{3k}k}{(2k+1)27^k}H_k^{(2)}\eq&-3(2-2p)B_{p^2-p-1}\l(\f13\r)\pmod{p^2}
\\\eq&-3B_{p-2}\l(\f13\r)\pmod p.\endalign$$
Below we show the first two congruences in (1.20) for $p>5$. (The case $p=5$ can be checked directly.)
Clearly $1/3\eq(2p^2+1)/3\pmod{p^2}$. Applying (1.14) and (1.15) with $a=1/3$ and $m=(2p^2+1)/3$, we obtain
$$\sum_{k=0}^{p-1}\f{\bi{2k}k\bi{3k}k}{27^k}H_k^{(2)}\eq-2\sum^{\f23(p^2-1)}\Sb k=1\\p\nmid k\endSb\f{(-1)^k}{k^2}\pmod{p^2}\tag3.11$$
and
$$-\f13\sum_{k=0}^{p-1}\f{\bi{2k}k\bi{3k}k}{(2k+1)27^k}H_k^{(2)}\eq-2\sum^{\f23(p^2-1)}\Sb k=1\\p\nmid k\endSb\f1{k^2}\pmod{p^2}.\tag3.12$$
On the other hand, (1.9) with $a=1/6$ and $m=(5p^2+1)/6$ yields
$$-\f23\sum_{k=0}^{p-1}\f{\bi{6k}{3k}\bi{3k}k}{(2k+1)432^k}H_k^{(2)}\eq-2\sum^{\f56(p^2-1)}\Sb k=1\\p\nmid k\endSb\f1{k^2}\pmod{p^2}.\tag3.13$$

Observe that
$$\align 2\sum^{\f23(p^2-1)}\Sb k=1\\p\nmid k\endSb\f{(-1)^k+1}{k^2}=&2\sum^{(p^2-1)/3}\Sb j=1\\p\nmid j\endSb\f2{(2j)^2}
\\\eq&\sum^{(p^2-1)/3}\Sb j=1\\p\nmid j\endSb\f1{(p^2-j)^2}=\sum^{p^2-1}\Sb k=1\\p\nmid k\endSb\f1{k^2}-\sum^{\f23(p^2-1)}\Sb k=1\\p\nmid k\endSb\f1{k^2}
\\\eq&-\sum^{\f23(p^2-1)}\Sb k=1\\p\nmid k\endSb\f1{k^2}\pmod{p^2}.
\endalign$$
Thus
$$2\sum^{\f23(p^2-1)}\Sb k=1\\p\nmid k\endSb\f{(-1)^k}{k^2}\eq-3\sum^{\f23(p^2-1)}\Sb k=1\\p\nmid k\endSb\f1{k^2}\pmod{p^2}\tag3.14$$
and
$$\sum^{(p^2-1)/3}\Sb k=1\\p\nmid k\endSb\f1{k^2}\eq-\sum^{\f23(p^2-1)}\Sb k=1\\p\nmid k\endSb\f1{k^2}\eq\f23\sum^{\f23(p^2-1)}\Sb k=1\\p\nmid k\endSb\f{(-1)^k}{k^2}\pmod{p^2}.\tag3.15$$
With the help of (3.8), we have
$$\align \sum^{(p^2-1)/3}\Sb k=1\\p\nmid k\endSb\f{(-1)^k}{k^2}=&\sum^{p^2-1}\Sb k=1\\p\nmid k\endSb\f{(-1)^k}{k^2}-\sum^{\f23(p^2-1)}\Sb k=1\\p\nmid k\endSb\f{(-1)^{p^2-k}}{(p^2-k)^2}
\\\eq&\sum^{p^2-1}\Sb k=1\\p\nmid k\endSb\f{(-1)^k+1}{k^2}+\sum^{\f23(p^2-1)}\Sb k=1\\p\nmid k\endSb\f{(-1)^k}{k^2}
\\\eq&\sum^{(p^2-1)/2}\Sb k=1\\p\nmid k\endSb\f{2}{(2k)^2}+\sum^{\f23(p^2-1)}\Sb k=1\\p\nmid k\endSb\f{(-1)^k}{k^2}\pmod{p^2}
\endalign$$
and hence
$$ \sum^{(p^2-1)/3}\Sb k=1\\p\nmid k\endSb\f{(-1)^k}{k^2}\eq\sum^{\f23(p^2-1)}\Sb k=1\\p\nmid k\endSb\f{(-1)^k}{k^2}\pmod{p^2}.\tag3.16$$
Adding (3.15) and (3.16) we get
$$\sum^{(p^2-1)/6}\Sb k=1\\p\nmid k\endSb\f2{(2k)^2}\eq\f53\sum^{\f23(p^2-1)}\Sb k=1\\p\nmid k\endSb\f{(-1)^k}{k^2}\pmod{p^2}$$
and hence
$$\align\sum^{\f56(p^2-1)}\Sb k=1\\p\nmid k\endSb\f1{k^2}=&\sum^{p^2-1}\Sb k=1\\p\nmid k\endSb\f1{k^2}-\sum^{(p^2-1)/6}\Sb k=1\\p\nmid k\endSb\f1{(p^2-k)^2}
\\\eq&-\sum^{(p^2-1)/6}\Sb k=1\\p\nmid k\endSb\f1{k^2}\eq-\f{10}3\sum^{\f23(p^2-1)}\Sb k=1\\p\nmid k\endSb\f{(-1)^k}{k^2}\pmod{p^2}.
\endalign$$
Combining this with (3.11)-(3.14) we obtain the first two congruences in (1.20). (When $p=5$, the second congruence in (1.20) can be verified directly.)
This concludes the proof. \qed

\heading{4. Some conjectural congruences}\endheading

In this section we pose some new conjectures on congruences, which are different from
the 100 conjectures in [Su19]. For any prime $p$, we use $\Z_p$ to denote
the ring of all $p$-adic integers.

\proclaim{Conjecture 4.1} Let $p$ be an odd prime. Then
$$\align\sum_{k=0}^{(p-1)/2}\f{(3k-1)4^k}{(2k-1)k^3\bi{2k}{k}^2}&\eq-\f74B_{p-3}\pmod{p},
\\ \sum_{k=1}^{(p-1)/2}\f{30k-11}{(2k-1)k^3\bi{2k}k^2}&\eq-8B_{p-3}\pmod p,
\endalign$$
and
$$\sum_{k=0}^{(p-1)/2}\f{(30k+11)\bi{2k}k^2}{k(2k+1)}
\eq8-8p-\f{112}3p^3B_{p-3}\pmod{p^4}$$
if $p>3$.
\endproclaim
\Remark\ 4.1. By Examples 72 and 11 of Chu and Zhang [CZ], we have
$$\sum_{k=1}^\infty\f{(3k-1)4^k}{(2k-1)k^3\bi{2k}k^2}=\f74\zeta(3)
\ \ \t{and}\ \ \sum_{k=1}^\infty\f{30k-11}{(2k-1)k^3\bi{2k}k^2}=4\zeta(3).$$

\proclaim{Conjecture 4.2} Let $p$ be an odd prime.

{\rm (i)} We have
$$\sum_{k=1}^{p-1}(-1)^k(7k+2)\bi{2k}k\f{\bi{3k}{k}}{2k+1}\eq7p^3B_{p-3}\pmod{p^4}$$
and
$$\sum_{k=1}^{(p-1)/2}(-1)^k(7k+2)\bi{2k}k\f{\bi{3k}{k}}{2k+1}
\eq3p\,q_p(2)-\f 32p^2q_p(2)^2
\pmod{p^3}.$$

{\rm (ii)} For any positive integer $n$, we have
$$\f1{(pn)^3}\sum_{k=n}^{pn-1}(-1)^k(7k+2)\bi{2k}k\f{\bi{3k}{k}}{2k+1}
\in\Z_p.$$
\endproclaim
\Remark\ 4.2. Chu and Zhang [CZ, Example 24] has the following equivalent form:
$$\sum_{k=1}^\infty\f{(-1)^k(7k-2)}{(2k-1)k^2\bi{2k}k\bi{3k}{k}}=-\f{\pi^2}{12}.$$

\proclaim{Conjecture 4.3}
{\rm (i)} Let $p$ be an odd prime. Then
$$\sum_{k=0}^{p-1}\f{(6k+1)\bi{2k}k\bi{4k}{2k}}{(2k+1)16^k}\eq\l(\f{-1}p\r)+3p^2E_{p-3}\pmod{p^3},$$
and
$$\sum_{k=0}^{(p-1)/2}\f{(6k+1)\bi{2k}k\bi{4k}{2k}}{(2k+1)16^k}\eq\l(\f{-1}p\r)\l(1+p\,q_p(2)-\f{p^2}2q_p(2)^2\r)
\pmod{p^3}$$
if $p>3$.

{\rm (ii)} For any odd prime $p$ and positive integer $n$, we have
$$\f1{(pn)^2}\(\sum_{k=0}^{pn-1}\f{(6k+1)\bi{2k}k\bi{4k}{2k}}{(2k+1)16^k}
-\l(\f{-1}p\r)\sum_{k=0}^{n-1}\f{(6k+1)\bi{2k}k\bi{4k}{2k}}{(2k+1)16^k}\)\in\Z_p.$$

{\rm (iii)} For any prime $p>3$, we have
$$\sum_{k=0}^{(p-1)/2}\f{(6k+1)\bi{2k}k\bi{4k}{2k}}{(2k+1)16^k}(4H_{2k}-H_k)
\eq\l(\f{-1}p\r)(4q_p(2)-p\,q_p(2)^2)\pmod{p^2}.$$
\endproclaim
\Remark\ 4.3. Chu and Zhang [CZ, Example 84] has the following equivalent form:
$$\sum_{k=1}^\infty\f{(6k-1)16^k}{(2k-1)k^2\bi{2k}k\bi{4k}{2k}}=8G,$$
where $G=\sum_{k=0}^\infty(-1)^k/(2k+1)^2$ is the Catalan constant.

\proclaim{Conjecture 4.4} Let $p$ be an odd prime.

{\rm (i)} We have
$$\sum_{k=1}^{(p-1)/2}\f{(5k+1)\bi{2k}k\bi{4k}{2k}}{(2k+1)(-16)^k}\eq
p\,q_p(2)-\f{p^2}2q_p(2)^2\pmod{p^3}$$
and
$$\sum_{k=1}^{p-1}\f{(5k+1)\bi{2k}k\bi{4k}{2k}}{(2k+1)(-16)^k}\eq\f 72p^3B_{p-3}\pmod{p^4}.$$

{\rm (ii)} For any positive integer $n$, we have
$$\f1{(pn)^3}\sum_{k=n}^{pn-1}\f{(5k+1)\bi{2k}k\bi{4k}{2k}}{(2k+1)(-16)^k}\in\Z_p.$$
\endproclaim
\Remark\ 4.4. Theorem 9 of Chu and Zhang [CZ] with $a=e=1$ and $b=c=d=1/2$ yields the identity
$$\sum_{k=1}^\infty\f{(5k-1)(-16)^k}{(2k-1)k^2\bi{2k}k\bi{4k}{2k}}=-\f{\pi^2}2.$$

\proclaim{Conjecture 4.5} Let $p$ be an odd prime.

{\rm (i)} If $p>3$, then
$$\sum_{k=0}^{p-1}\f{(22k+1)\bi{3k}k\bi{6k}{3k}}{(2k+1)256^k}\eq\l(\f{-1}p\r)+15p^2E_{p-3}\pmod{p^3}$$
and
$$\sum_{k=0}^{(p-1)/2}\f{(22k+1)\bi{3k}k\bi{6k}{3k}}{(2k+1)256^k}\eq\l(\f{-1}p\r)
\l(1+\f3{80}p^5B_{p-3}\r)\pmod{p^6}.$$

{\rm (ii)} For any positive integer $n$, we have
$$\f1{(pn)^2}\(\sum_{k=0}^{pn-1}\f{(22k+1)\bi{3k}k\bi{6k}{3k}}{(2k+1)256^k}
-\l(\f{-1}p\r)\sum_{k=0}^{n-1}\f{(22k+1)\bi{3k}k\bi{6k}{3k}}{(2k+1)256^k}\)\in\Z_p.$$
If $p\not=5$, then for any positive odd integer $n$ we have
$$\f1{(pn)^5}\(\sum_{k=0}^{(pn-1)/2}\f{(22k+1)\bi{3k}k\bi{6k}{3k}}{(2k+1)256^k}
-\l(\f{-1}p\r)\sum_{k=0}^{(n-1)/2}\f{(22k+1)\bi{3k}k\bi{6k}{3k}}{(2k+1)256^k}\)\in\Z_p.$$
\endproclaim
\Remark\ 4.5. Chu and Zhang [CZ, Example 50] proved that
$$\sum_{k=1}^\infty\f{(22k-1)256^k}{(2k-1)k^2\bi{3k}k\bi{6k}{3k}}=128G.$$

\proclaim{Conjecture 4.6} Let $p$ be an odd prime.

{\rm (i)} We have
$$\sum_{k=0}^{(p-1)/2}\frac{(10k-1)\binom{3k}k\binom{6k}{3k}}{(2k+1)512^k}
\equiv\left(\frac{-2}p\right)\bigg(\frac98p^2q_p(2)^2-\frac32p\,q_p(2)-1\bigg)\pmod{p^3},$$
and
$$\sum_{k=0}^{p-1}\frac{(10k-1)\binom{3k}k\binom{6k}{3k}}{(2k+1)512^k}
\equiv-\left(\frac{-2}p\right)+\frac{15}{16}p^2E_{p-3}\left(\frac14\right)
\pmod{p^3}.$$

{\rm (ii)} For any positive integer $n$, we have $$\frac1{(pn)^2}\bigg(\sum_{k=0}^{pn-1}\frac{(10k-1)\binom{3k}k\binom{6k}{3k}}{(2k+1)512^k}-\left(\frac{-2}p
\right)\sum_{k=0}^{n-1}\frac{(10k-1)\binom{3k}k\binom{6k}{3k}}{(2k+1)512^k}\bigg)\in\Z_p.$$
\endproclaim
\Remark\ 4.6. We also conjecture that
$$\sum_{k=0}^\infty\frac{(10k-1)\binom{3k}k\binom{6k}{3k}}{(2k+1)512^k}=0$$
and
$$\sum_{k=0}^\infty\f{(10k-1)\bi{3k}k\bi{6k}{3k}}{(2k+1)512^k}\l(H_{2k}-\f23H_k\r)=\sum_{k=0}^\infty
\f{\bi{3k}k\bi{6k}{3k}}{(2k+1)^2512^k}.$$

\proclaim{Conjecture 4.7} Let $p$ be any odd prime.

{\rm (i)} We have
$$\sum_{k=0}^{(p-1)/2}\f{(7k+1)\bi{3k}k\bi{6k}{3k}}{(2k+1)(-4)^k\bi{2k}k}
\eq\l(\f{-1}p\r)\l(1+\f32p\,q_p(2)\r)\pmod{p^3}$$
and
$$\sum_{k=0}^{p-1}\f{(7k+1)\bi{3k}k\bi{6k}{3k}}{(2k+1)(-4)^k\bi{2k}k}
\eq\l(\f{-1}p\r)-15p^2E_{p-3}\pmod{p^3}.$$

{\rm (ii)} For any positive integer $n$, we have
$$\f1{(pn)^2}\(\sum_{k=0}^{pn-1}\f{(7k+1)\bi{3k}k\bi{6k}{3k}}{(2k+1)(-4)^k\bi{2k}k}
-\l(\f{-1}p\r)\sum_{k=0}^{n-1}\f{(7k+1)\bi{3k}k\bi{6k}{3k}}{(2k+1)(-4)^k\bi{2k}k}\)\in\Z_p.$$
\endproclaim
\Remark\ 4.7. Chu and Zhang [CZ, Example 27] has the following equivalent form:
$$\sum_{k=1}^\infty\f{(7k-1)(-4)^k\bi{2k}k}{(2k-1)k\bi{3k}k\bi{6k}{3k}}=-\f{\pi}4.$$

\proclaim{Conjecture 4.8} Let $p>3$ be a prime.

{\rm (i)} We have
$$\sum_{k=0}^{(p-1)/2}\f{\bi{3k}k\bi{6k}{3k}}{216^k\bi{2k}k}
\eq\l(\f 6p\r)\l(1+\f p6q_p(2)\r)\pmod{p^2}.$$

{\rm (ii)} For any positive integer $n$, we have
$$\f1{pn}\(\sum_{k=0}^{pn-1}\f{\bi{3k}k\bi{6k}{3k}}{216^k\bi{2k}k}
-\l(\f{6}p\r)\sum_{k=0}^{n-1}\f{\bi{3k}k\bi{6k}{3k}}{216^k\bi{2k}k}\)\in\Z_p$$
and
$$\f1{pn}\(\sum_{k=0}^{pn-1}\f{\bi{3k}k\bi{6k}{3k}}{(2k+1)216^k\bi{2k}k}
-\l(\f{2}p\r)\sum_{k=0}^{n-1}\f{\bi{3k}k\bi{6k}{3k}}{(2k+1)216^k\bi{2k}k}\)\in\Z_p.$$
\endproclaim
\Remark\ 4.8. By [OLBC, 15.4.30], we have
$$\sum_{k=0}^\infty\f{\bi{3k}k\bi{6k}{3k}}{216^k\bi{2k}k}=\f{\sqrt6}2
\ \ \t{and}\ \ \sum_{k=0}^\infty\f{\bi{3k}k\bi{6k}{3k}}{(2k+1)216^k\bi{2k}k}=\f{3}4\sqrt2.$$
The author guessed the identities
$$\sum_{k=0}^\infty\f{\bi{3k}k\bi{6k}{3k}}{216^k\bi{2k}k}\sum_{0\ls j<k}\f1{2j+1}=\f{3\sqrt6}8\log\f43$$
and $$\sum_{k=0}^\infty\f{\bi{3k}k\bi{6k}{3k}}{(2k+1)216^k\bi{2k}k}\sum_{j=0}^k\f1{2j+1}
=\f{3\sqrt2}{16}\l(6+\log\f4{27}\r),$$
which were later confirmed by his PhD student W. Xia.

\proclaim{Conjecture 4.9} Let $p$ be an odd prime.

{\rm (i)} We have
$$\sum_{k=0}^{(p-1)/2}\f{\bi{3k}k\bi{6k}{3k}}{(2k+1)8^k\bi{2k}k}\eq\l(\f2p\r)\l(\f{1+(\f{-1}p)}2-\f34p\,q_p(2)\r)
\pmod{p^2},$$
$$\sum_{k=0}^{p-1}\f{(5k+1)\bi{3k}k\bi{6k}{3k}}{(2k+1)8^k\bi{2k}k}
\eq\l(\f {-2}p\r)-\f{15}{16}p^2E_{p-3}\l(\f14\r)\pmod{p^3}$$
and
$$\sum_{k=0}^{(p-1)/2}\f{(5k+1)\bi{3k}k\bi{6k}{3k}}{(2k+1)8^k\bi{2k}k}
\eq\l(\f {2}p\r)\l(\f{3+(\f{-1}p)}4+\f38 p\,q_p(2)\r)\pmod{p^2}.$$

{\rm (ii)} For any positive integer $n$, we have
$$\f1{(pn)^2}\(\sum_{k=0}^{pn-1}\f{(5k+1)\bi{3k}k\bi{6k}{3k}}{(2k+1)8^k\bi{2k}k}
-\l(\f{-2}p\r)\sum_{k=0}^{n-1}\f{(5k+1)\bi{3k}k\bi{6k}{3k}}{(2k+1)8^k\bi{2k}k}\)\in\Z_p.$$
\endproclaim
\Remark\ 4.9. We note that
$$\sum_{k=0}^{p-1}\f{\bi{3k}k\bi{6k}{3k}}{(2k+1)8^k\bi{2k}k}\eq\l(\f 2p\r)\l(2\l(\f{-1}p\r)-1\r)\pmod p$$
for any odd prime $p$.

\proclaim{Conjecture 4.10} Let $p$ be an odd prime.

{\rm (i)} We have
$$\sum_{k=1}^{(p-1)/2}\f{(74k+7)\bi{2k}k\bi{3k}k\bi{6k}{3k}}{(2k+1)4096^k}
\eq-\f 94pH_{p-1}\pmod{p^5}$$
if $p\not=5$.
and
$$\sum_{k=1}^{p-1}\f{(74k+7)\bi{2k}k\bi{3k}k\bi{6k}{3k}}{(2k+1)4096^k}
\eq -\f{35}2p^3B_{p-3}\pmod{p^4}.$$

{\rm (ii)} For any positive integer $n$, we have
$$\f1{(pn)^3}\sum_{k=n}^{pn-1}\f{(74k+7)\bi{2k}k\bi{3k}k\bi{6k}{3k}}{(2k+1)4096^k}\in\Z_p$$
and
$$\f1{(p(2n-1))^3}\sum_{k=n}^{pn-(p+1)/2}\f{(74k+7)\bi{2k}k\bi{3k}k\bi{6k}{3k}}{(2k+1)4096^k}\in\Z_p.$$
\endproclaim
\Remark\ 4.10. Chu and Zhang [CZ, Example 60] has the following equivalent form:
$$\sum_{k=0}^{\infty}\f{(74k+7)\bi{2k}k\bi{3k}k\bi{6k}{3k}}{(2k+1)4096^k}=8.$$

\proclaim{Conjecture 4.11} Let $p>3$ be a prime.

{\rm (i)} We have
$$\sum_{k=0}^{(p-1)/2}(42k^2+27k+4)\f{\bi{2k}k\bi{3k}k\bi{4k}{2k}}{(2k+1)(-9)^k}
\eq2p\l(1+\l(\f p3\r)\r)-2p^2q_p(3)\pmod{p^3}$$
and
$$\sum_{k=0}^{p-1}(42k^2+27k+4)\f{\bi{2k}k\bi{3k}k\bi{4k}{2k}}{(2k+1)(-9)^k}
\eq 4p\l(\f p3\r)+\f 83 p^3B_{p-2}\l(\f13\r)\pmod{p^4}.$$

{\rm (ii)} For any positive integer $n$, we have
$$\f1{(pn)^3}\(\sum_{k=0}^{pn-1}\f{P(k)\bi{2k}k\bi{3k}k\bi{4k}{2k}}{(2k+1)(-9)^k}
-p\l(\f p3\r)\sum_{k=0}^{n-1}\f{P(k)\bi{2k}k\bi{3k}k\bi{4k}{2k}}{(2k+1)(-9)^k}\)\in\Z_p,$$
where $P(k)=42k^2+27k+4$.
\endproclaim
\Remark\ 4.11. We also conjecture that
$$\sum_{k=1}^\infty\frac{(42k^2-27k+4)(-9)^k}{k^3(2k-1)\binom{2k}k\binom{3k}k\binom{4k}{2k}}=-6
\sum_{k=1}^\infty\frac{(\frac k3)}{k^2}$$
and
$$\sum_{k=1}^\infty\frac{(-9)^k\left((42k^2-27k+4)(H_{3k-1}-H_{k-1})-10k+3\right)}{k^3(2k-1)
\binom{2k}k\binom{3k}k\binom{4k}{2k}}=-\frac{8\pi^3}{27\sqrt3}.$$

\proclaim{Conjecture 4.12} Let $p>3$ be a prime.

{\rm (i)} We have
$$\sum_{k=0}^{(p-1)/2}\f{(35k^2+29k+6)\bi{4k}k}{(3k+1)(3k+2)3^k}\eq\l(\f p3\r)\pmod p$$
and
$$\sum_{k=0}^{p-1}\f{(35k^2+29k+6)\bi{4k}k}{(3k+1)(3k+2)3^k}\eq3\l(\f p3\r)-\f 83p^2B_{p-2}\l(\f13\r)
\pmod{p^2}.$$

{\rm (ii)} For any positive integer $n$, we have
$$\f1{(pn)^2}\(\sum_{k=0}^{pn-1}\f{(35k^2+29k+6)\bi{4k}k}{(3k+1)(3k+2)3^k}-\l(\f p3\r)\sum_{k=0}^{n-1}\f{(35k^2+29k+6)\bi{4k}k}{(3k+1)(3k+2)3^k}\)\in\Z_p.$$
\endproclaim
\Remark\ 4.12. We also conjecture that
$$\sum_{k=1}^\infty\f{(35k^2-29k+6)3^k}{(3k-1)(3k-2)k\bi{4k}k}=\sqrt3\,\pi.$$

\proclaim{Conjecture 4.13}
{\rm (i)} Let $p$ be an odd prime. Then
$$\sum_{k=0}^{(p-3)/2}(145k^2+104k+18)\f{\bi{2k}k\bi{3k}k^2}{2k+1}\eq 9p+28p^2\l(\f{-1}p\r)\pmod{p^3},$$
and
$$\sum_{k=0}^{p-1}(145k^2+104k+18)\f{\bi{2k}k\bi{3k}k^2}{2k+1}
\eq18p-288p^2H_{p-1}+\f{4572}5p^6B_{p-5}\pmod{p^7}$$
if $p>3$.

{\rm (ii)} For any prime $p$ and positive integer $n$, we have
$$\f1{(pn)^2}\(\sum_{k=0}^{pn-1}(145k^2+104k+18)\f{\bi{2k}k\bi{3k}k^2}{2k+1}
-p\sum_{k=0}^{n-1}(145k^2+104k+18)\f{\bi{2k}k\bi{3k}k^2}{2k+1}\)\in\Z_p.$$
\endproclaim
\Remark\ 4.13. We also conjecture that
$$\sum_{k=1}^\infty\frac{145k^2-104k+18}{k^3(2k-1)\binom{2k}k\binom{3k}k^2}=\frac{\pi^2}3.$$

\proclaim{Conjecture 4.14} Let $p$ be an odd prime.

{\rm (i)} We have
$$\sum_{k=0}^{p-1}(42k^2+23k+3)\f{\bi{2k}k\bi{4k}{2k}^2}{(2k+1)16^k}\eq3p+\f{63}2p^4B_{p-3}\pmod{p^5},$$
and
$$\sum_{k=0}^{(p-1)/2}(42k^2+23k+3)\f{\bi{2k}k\bi{4k}{2k}^2}{(2k+1)16^k}\eq p\l(1+2\l(\f{-1}p\r)\r)-2\l(\f{-1}p\r)p^3E_{p-3}\pmod{p^4}$$
if $p>3$.

{\rm (ii)} For any positive integer $n$, we have
$$\f1{(pn)^4}\(\sum_{k=0}^{pn-1}(42k^2+23k+3)\f{\bi{2k}k\bi{4k}{2k}^2}{(2k+1)16^k}
-p\sum_{k=0}^{n-1}(42k^2+23k+3)\f{\bi{2k}k\bi{4k}{2k}^2}{(2k+1)16^k}\)\in\Z_p.$$
\endproclaim
\Remark\ 4.14. We also conjecture that
$$\sum_{k=1}^\infty\frac{(42k^2-23k+3)16^k}{k^3(2k-1)\binom{2k}k\binom{4k}{2k}^2}=\frac{\pi^2}2.$$

\proclaim{Conjecture 4.15} Let $p$ be an odd prime.

{\rm (i)} We have
$$\sum_{k=0}^{(p-3)/2}(92k^2+61k+9)\f{\bi{2k}k\bi{3k}k\bi{4k}{2k}}{(2k+1)64^k}
\eq 6p+16p^2\l(\f{-1}p\r)\pmod{p^3}$$
and
$$\sum_{k=0}^{p-1}(92k^2+61k+9)\f{\bi{2k}k\bi{3k}k\bi{4k}{2k}}{(2k+1)64^k}
\eq 9p+\f{63}2p^4B_{p-3}\pmod{p^3}.$$

{\rm (ii)} For any positive integer $n$, we have
$$\f1{(pn)^4}\(\sum_{k=0}^{pn-1}(92k^2+61k+9)\f{\bi{2k}k\bi{3k}k\bi{4k}{2k}}{(2k+1)64^k}
-p\sum_{k=0}^{n-1}(92k^2+61k+9)\f{\bi{2k}k\bi{3k}k\bi{4k}{2k}}{(2k+1)64^k}\)\in\Z_p.$$
\endproclaim
\Remark\ 4.15. We also conjecture that
$$\sum_{k=1}^\infty\frac{(92k^2-61k+9)64^k}{k^3(2k-1)\binom{2k}k\binom{3k}k\binom{4k}{2k}}=8\pi^2.$$

\proclaim{Conjecture 4.16}
Let $p$ be an odd prime.

{\rm (i)} We have
$$\sum_{k=0}^{(p-3)/2}\f{(592k^3+580k^2+112k-3)\bi{3k}k\bi{6k}{3k}}{(2k+1)(4k+1)(4k+3)1024^k}
\eq\l(\f{-1}p\r)(8p^2-4)
\pmod{p^3}$$
and
$$\sum_{k=0}^{p-1}\f{(592k^3+580k^2+112k-3)\bi{3k}k\bi{6k}{3k}}{(2k+1)(4k+1)(4k+3)1024^k}\eq15p^2E_{p-3}-\l(\f{-1}p\r)\pmod{p^3}.$$

{\rm (ii)} For any positive integer $n$, we have
$$\f1{(pn)^2}\(\sum_{k=0}^{pn-1}f(k)
-\l(\f{-1}p\r)\sum_{k=0}^{n-1}f(k)\)
\in\Z_p,$$
where
$$f(k)=\f{(592k^3+580k^2+112k-3)\bi{3k}k\bi{6k}{3k}}{(2k+1)(4k+1)(4k+3)1024^k}.$$
\endproclaim
\Remark\ 4.16. We also conjecture that
$$\sum_{k=0}^\infty\f{P(k)\bi{3k}k\bi{6k}{3k}}{(2k+1)(4k+1)(4k+3)1024^k}=0$$
and
$$\sum_{k=0}^\infty\f{\bi{3k}k\bi{6k}{3k}(P(k)(2H_{6k}-H_{3k}-3H_{2k}+H_k)
-Q(k)/(2k+1))}{(2k+1)(4k+1)(4k+3)1024^k}=0,$$
where
$$P(k)=592k^3+580k^2+112k-3\ \ \t{and}\ \ Q(k)=2096k^3+2076k^2+400k-5.$$

Actually, we also have many other similar conjectures.

\Ack. The author thanks his PhD student Wei Xia for helpful comments.

\enddocument